
\documentstyle{amsppt}
\NoBlackBoxes
\mag=\magstep1 
\font \smallfont = cmr8 
\hsize=6.5 true in
\vsize=9.0 true in
\voffset=-0.5 true in
\define\docbaseline{\baselineskip=19 pt} 

\define\reals{\Bbb R}
\define\setof#1#2{\{\,#1\mid#2\,\}}
\define\set#1{\{#1\}}
\define\bigsetof#1#2{\bigl\{\,#1 \bigm| #2\,\bigr\}}
\define\gl{\operatorname{gl}}
\define\Sym{\operatorname{Sym}}
\define\Pan{\operatorname{Pan}}
\define\dottimes{\, \dot{\times}\, }
\define\modulus{\operatorname{mod}}
\define\support{\operatorname{supp}}
\define\mymod#1{\ (\modulus \, #1)}


\topmatter

\title{Birkhoff's Theorem for Panstochastic Matrices}\endtitle

\author{Dean~Alvis and Michael~Kinyon}\endauthor

\endtopmatter

\document
\headline={\hss \tenrm --\ \folio\ -- \hss}
\footline={\hfill}

\docbaseline


\bigskip
\proclaim{1. Introduction}
\endproclaim

An $n \times n$ matrix with nonnegative real entries is 
{\it doubly stochastic} if the sum of the entries along
any of its rows or columns is equal to $1$. A doubly
stochastic matrix is {\it panstochastic} if 
the sum of the entries along any downward diagonal or
upward diagonal, either broken or unbroken, is equal to
$1$. For example, the matrices
$$
{1 \over 5}
\left [
\matrix
1 & 1 & 1 & 1 & 1 \cr
1 & 1 & 1 & 1 & 1 \cr
1 & 1 & 1 & 1 & 1 \cr
1 & 1 & 1 & 1 & 1 \cr
1 & 1 & 1 & 1 & 1 \cr
\endmatrix
\right ]
,
\qquad
\left [
\matrix
1 & 0 & 0 & 0 & 0 \cr
0 & 0 & 0 & 1 & 0 \cr
0 & 1 & 0 & 0 & 0 \cr
0 & 0 & 0 & 0 & 1 \cr
0 & 0 & 1 & 0 & 0 \cr
\endmatrix
\right ],
\qquad
{1 \over 60}
\left [
\matrix
1 & 13 & 20 & 7 & 19 \cr
22 & 9 & 16 & 3 & 10 \cr
18 & 0 & 12 & 24 & 6 \cr
14 & 21 & 8 & 15 & 2 \cr
5 & 17 & 4 & 11 & 23 \cr
\endmatrix
\right ]
$$
are panstochastic. 

A linear combination 
is called {\it convex} if 
the coefficients are nonnegative and their sum 
is equal to $1$.
In \cite{1}, Birkhoff showed that every 
doubly stochastic matrix can be expressed as a convex 
combination of permutation matrices. 
Related results for integral matrices had been 
obtained earlier by K\H{o}nig \cite{8} and 
Egerv\'ary \cite{5}. Birkhoff's 
theorem has been generalized in various ways; for 
example, Schneider obtained the result for matrices 
with entries in lattice-ordered abelian groups 
\cite{11}.

Does the analogue of
Birkhoff's theorem hold for panstochastic matrices? 
Our first main result is that this is the case when $n=5$. 

\proclaim{Theorem~1.1}
A $5 \times 5$ real matrix 
is panstochastic if and only if it 
is a convex combination of panstochastic permutation matrices. 
\endproclaim

Because no permutation matrix is a convex 
combination of other permutation matrices, the $n \times n$
doubly-stochastic matrices form a generalized polyhedron whose
vertices are the $n \times n$ permutation matrices
\cite{9}. {Theorem~1.1}
has the following geometric interpretation.

\proclaim{Corollary~1.1}
The set of all real panstochastic $5 \times 5$ matrices 
forms a convex polyhedron whose vertices 
are panstochastic permutation
matrices.
\endproclaim

The analogue of Birkhoff's theorem
for panstochastic matrices clearly holds if $n = 1$. 
Our second main result shows that {Theorem~1.1} does not
generalize for $n>1$, $n \ne 5$.

\proclaim{Theorem~1.2}
If $n > 1$ and $n \ne 5$, then there is some $n \times n$ panstochastic
matrix that is not a convex combination of panstochastic permutation
matrices.
\endproclaim

We adopt the following notation and terminology: 
An $n \times n$ matrix $A$ over a field is {\it panmagic} if
the sums of the entries along all of the rows, columns, downward
diagonals, and upward diagonals, both broken and unbroken, of $A$ are
equal, in which case the common value of the sums is called the
{\it magic number} of $A$, denoted $\mu(A)$. The set of all
$n \times n$ panmagic matrices over $F$, denoted $\Pan(n,F)$, is a
subspace of $\gl(n,F)$, the vector space of all $n \times n$ matrices
over $F$. We index the rows and columns of an $n \times n$ matrix 
by the elements of $\Omega_n = \set{0,1,\dots,n-1}$. 
We index the diagonals 
in the following way: for $k \in \Omega_n$, the $k$th
upward diagonal contains the $(i,j)$ entry whenever
$i+j \equiv k \ (\modulus \, n)$,
and the $k$th downward diagonal contains the $(i,j)$ entry whenever
$i-j \equiv k \ (\modulus \, n)$.


\bigskip
\proclaim{2. Panmagic Permutations}
\endproclaim

Let $\Sym(\Omega_n)$ be the group of all permutations on $\Omega_n$,
and let $F$ be a field. For $\pi \in \Sym(\Omega_n)$, let
$P_\pi$ denote the corresponding permutation matrix: the $(i,j)$
entry of $P_\pi$ is $1$ if $i = \pi(j)$ and $0$ otherwise. We say
$\pi \in \Sym(\Omega_n)$ is a {\it panmagic permutation} of degree $n$
if $P_\pi \in \Pan(n,F)$. Denote by $\Pi_n$ the set of all panmagic
permutations of degree $n$. 
These definitions are independent of $F$; see {Lemma~2.1}.
In particular, if $F=\reals$, then $\pi \in \Sym(\Omega_n)$ is panmagic 
if and only if $P_\pi$ is panstochastic.

If $\pi \in \Sym(\Omega_n)$, then the sum of the entries along
any row or column of $P_\pi$ is equal to $1$. The sums of the entries 
along the diagonals of $P_\pi$ are all equal to $1$ if and only if the 
congruences
$\pi(j) - j \equiv k \mymod{n}$
and
$\pi(j) + j \equiv k \mymod{n}$
have unique solutions for $j \in \Omega_n$ whenever $k \in \Omega_n$.
We record this observation as 

\proclaim{Lemma~2.1}
Let $\pi \in \Sym(\Omega_n)$. Then $\pi$ is panmagic if and only if
there exist $\lambda$, $\rho \in \Sym(\Omega_n)$ such that
$$
\pi(j)-j \equiv \lambda(j) \ (\modulus \ n)
\qquad
\text{and}
\qquad
\pi(j)+j \equiv \rho(j) \ (\modulus \ n)
$$
for all $j \in \Omega_n$.
\endproclaim


Let $a$ be an integer that is relatively prime to $n$, and let $b$ be an
arbitrary integer. Denote by $\pi_{ax+b}$ the permutation of $\Omega_n$
that sends $j$ to $a j+b \ (\modulus\ n)$. Such a permutation is 
called an {\it affine permutation} of $\Omega_n$. If $\pi = \pi_{ax+b}$
is an affine permutation, then we sometimes denote 
$P_{\pi}$ by 
$P_{ax+b}$. Let $\Lambda_n $ be the set of all affine panmagic
permutations of degree $n$. The following is an immediate consequence
of {Lemma~2.1}.


\proclaim{Lemma~2.2}
Suppose $a$ is relatively prime to $n$. Then the affine permutation
$\pi_{ax + b}$ is panmagic if and only if $a+1$ and $a-1$ are also
relatively prime to $n$.
\endproclaim


Squaring the congruences in {Lemma~2.1} and summing over $j$ gives the
following result; see \cite{3}, 
\cite{6}, and \cite{7}.

\proclaim{Theorem~2.1}
There is some panmagic permutation of degree $n$ if and only if
$\gcd(n,6)=1$.
\endproclaim


This result is already enough to prove one case of {Theorem~1.2}.
If $\gcd(n,6)>1$, then the panstochastic matrix with all entries 
equal to $1/n$ is not a convex combination of panmagic permutation
matrices. We can therefore rule out this case.


We next determine which $n$ satisfy $\Pi_n = \Lambda_n$. Let $n=5$.
If $\pi \in \Pi_5$ and $\pi(0)=b$, then the $b$th entry of the zeroth
column of $P_{\pi}$ must be $1$. This forces all other entries in the
zeroth column, the $b$th row, the $b$th upward diagonal, and the
$(n-b)$th downward diagonal to be $0$. Consideration of the remaining
eight entries shows that either $\pi = \pi_{2x+b}$ or
$\pi = \pi_{3x+b}$. Therefore every panmagic permutation of degree $5$
is affine: $\Pi_5 = \Lambda_5$. Similar reasoning shows that 
$\Pi_7 = \Lambda_7$ and $\Pi_{11} = \Lambda_{11}$. 
On the other hand, 
$$
\Pi_n \ne \Lambda_n
\quad
\text{whenever}
\quad
n \ge 13
\quad
\text{and}
\quad
\gcd(n,6)=1.
$$
Indeed, Bruen and Dixon 
have constructed non-affine panmagic
permutations of degree $n$ 
whenever $n$ is prime and $n \ge 13$ \cite{2}.
See also \cite{7},
in which non-affine panmagic permutations are constructed
whenever $\gcd(n,6)=1$, $n$ is composite, and $n$ is not square-free.
We construct a variation of Bruen and Dixon's
example that applies when $n$ is composite
and $\gcd(n,6)=1$, as follows.
Let $p$ be a prime divisor of $n$. Define
$S = \setof{k \in \Omega_n}{k \equiv 0 \mymod{p}}$ 
and 
$T = \Omega_n \setminus S$. 
The permutation
$\pi : \Omega_n \rightarrow \Omega_n$ defined by
$$
\pi ( x ) =
\cases
2 x & \text{if $x \in S$,} \cr
3 x & \text{if $x \in T$,} \cr
\endcases
$$
is then panmagic and non-affine.


\bigskip
\proclaim{3. Kronecker and Wreath Products}
\endproclaim

In this section we consider methods for constructing new panmagic
matrices from existing ones. Assume that $m$ and $n$ are positive
integers and $F$ is a field. Suppose
$A_0, \dots, A_{n-1}\in \gl(m,F)$ and $B \in \gl(n,F)$. Define a
matrix $(A_0,\dots,A_{n-1}) \wr B \in \gl( m n, F)$ as follows: if 
$A_s = \left [ a^s_{i,j} \right ]_{i,j \in \Omega_m}$ 
and
$B = \left [ b_{r,s} \right]_{r,s \in \Omega_n}$, 
then the $(i n +r, j n +s)$ entry of $(A_0,\dots,A_{n-1}) \wr B$
is
$$
\left ( (A_0,\dots,A_{n-1}) \wr B \right )_{i n + r, j n + s}
=
a^s_{i,j} \, b_{r,s}
$$
for $i, j \in \Omega_m$, $r, s \in \Omega_n$. Notice that if
$A \in \gl(m,F)$, then
$ (A,\dots,A) \wr B = A \dottimes B,$
where $A \dottimes B$ is the Kronecker product of $A$ and $B$
\cite{4, p. 250}. The next result can be deduced from
\cite{10, Theorem~2.7 }.


\proclaim{Theorem~3.1}
Suppose $A_0, \dots A_{n-1}, \in \Pan(m,F)$ satisfy
$\mu(A_0)= \cdots = \mu(A_{n-1}),$
and let $B \in \Pan(n,F)$. Then 
$(A_0,\dots,A_{n-1}) \wr B \in \Pan(m n, F)$ and 
$\mu \left ( (A_0,\dots,A_{n-1}) \wr B \right )
= \mu(A_0) \mu(B)$.
\endproclaim


\proclaim{Corollary~3.1}
If $A \in \Pan(m,F)$ and $B \in \Pan(n,F)$,
then
$A \dottimes B \in \Pan(m n,F)$ and
$\mu(A \dottimes B) = 
\mu(A) \mu(B)$.
\endproclaim


Observe that if $\lambda_0, \dots, \lambda_{n-1} \in \Sym(\Omega_m)$
and $\rho \in \Sym(\Omega_n)$, then the matrix
$P_\pi = (P_{\lambda_0},\dots,P_{\lambda_{n-1}}) \wr P_\rho$
is a permutation matrix: the underlying permutation
$\pi \in \Sym(\Omega_{m n})$ is given by 
$$
\pi ( j n + s ) = \lambda_s (j) n + \rho (s)
\quad
\text{for}
\ j \in \Omega_m, s \in \Omega_n.
$$
Denote this permutation $\pi$ by
$(\lambda_0,\dots,\lambda_{n-1}) \wr \rho$. The collection of all
such permutations is a subgroup of $\Sym(\Omega_{m n})$ 
that is isomorphic to
the wreath product $\Sym(\Omega_m) \wr \Sym(\Omega_n)$.


\proclaim{Theorem~3.2}
If $\lambda_0, \dots, \lambda_{n-1} \in \Sym(\Omega_m)$ and
$\rho \in \Sym(\Omega_n)$,
then
$
(\lambda_0,\dots,\lambda_{n-1})~\wr~\rho \in \Pi_{m n}
$
if and only if $\lambda_0, \dots, \lambda_{n-1} \in \Pi_m$ and
$\rho \in \Pi_n$. 
\endproclaim

\demo{Proof}
Put $\pi = (\lambda_0,\dots,\lambda_{n-1}) \wr \rho$. If
$\lambda_0, \dots, \lambda_{n-1} \in \Pi_m$ and $\rho \in \Pi_n$, 
then $P_\pi = (P_{\lambda_0},\dots,P_{\lambda_{n-1}}) \wr P_\rho 
\in \Pan(m n, F)$ 
by {Theorem~3.1}, 
and thus $\pi \in \Pi_{m n}$.
Conversely, suppose $\rho \not \in \Pi_n$. By {Lemma~2.1} there is
some $\delta \in \left \{ -1,1 \right \}$ such that the mapping
$s \mapsto \rho(s) + \delta s \ (\modulus \, n)$ is not surjective from
$\Omega_n$ to $\Omega_n$. Hence the mapping
$x \mapsto \pi(x) + \delta x \ (\modulus \, m n)$ is not surjective from
$\Omega_{m n}$ to $\Omega_{m n}$ because
$\pi(j n + s) + \delta (j n + s) \equiv \rho(s) + \delta s
\ (\modulus \, n)$, and so $\pi \not \in \Pi_{m n}$. On the other
hand, suppose $\lambda_s \not \in \Pi_m$ for some $s \in \Omega_n$,
so $j \mapsto \lambda_s(j) + \delta j \ (\modulus \, m)$ is not
injective for some $\delta \in \left \{ -1,1 \right \}$. For this $s$
and $\delta$, the mapping
$j \mapsto \pi(j n + s) + \delta (j n + s) \ (\modulus \, m n)$
is not injective from $\Omega_m$ to $\Omega_{m n}$, and hence the
mapping $x \mapsto \pi(x) + \delta x \ (\modulus \, m n )$ is not
injective from $\Omega_{m n}$ to $\Omega_{m n}$. Therefore
$\pi \not \in \Pi_{m n}$.
\qed\enddemo


For $\lambda \in \Sym(\Omega_m)$, $\rho \in \Sym(\Omega_n)$, define
$\lambda \dottimes \rho =
(\lambda, \dots, \lambda) \wr \rho \in \Sym(\Omega_{m n})$.
Thus
$$
(\lambda \dottimes \rho) ( j n + s ) = \lambda(j) n + \rho(s)
$$
for $j \in \Omega_m$, $s \in \Omega_n$. The following result is a
special case of {Theorem~3.2}.

\proclaim{Corollary~3.2}
If $\lambda \in \Sym(\Omega_m)$ and $\rho \in \Sym(\Omega_n)$ are
panmagic, then $\lambda \dottimes \rho$ is panmagic.
\endproclaim


For $A=\left[ a_{i,j} \right]_{i,j\in\Omega_n} \in \gl(n,F)$, 
define 
$
\support A=
\bigsetof {(i,j) \in \Omega_n \times \Omega_n} {a_{i,j} \neq 0}
$.
The next result is used in the proof of {Lemma~4.2}.

\proclaim{Theorem~3.3}
Suppose $A \in \gl(m,F)$, $\rho \in \Sym(\Omega_n)$, 
$\pi \in \Sym(\Omega_{m n})$, and
$\support P_\pi \subseteq \support (A \dottimes P_\rho)$.
Then 
$
\pi = (\lambda_0, \dots, \lambda_{n-1}) \wr \rho
$
for some $\lambda_0, \dots, \lambda_{n-1} \in \Sym(\Omega_m)$,
and
$
\support P_{\lambda_s} \subseteq \support A
\text{ for all }
s \in \Omega_n.
$
\endproclaim

\demo{Proof}
Let $A = \left [ a_{i,j} \right ]_{i,j \in \Omega_m}$. Suppose
$i,j \in \Omega_m$, $r, s \in \Omega_n$ and $\pi(j n + s) = i n + r$. 
Thus the $(i n + r, j n + s)$ entry of $P_\pi$ is equal to $1$, so
the $(i n + r, j n + s)$ entry of $A \dottimes P_\rho$ is nonzero.
Hence $a_{i,j} (P_\rho)_{r,s} \ne 0$, and so $\rho(s) = r$. Thus there
are unique mappings $\lambda_0, \dots, \lambda_{n-1}$ from $\Omega_m$
to $\Omega_m$ such that
$
\pi( j n + s ) = \lambda_s(j) n + \rho( s ).
$
Since $\pi$ is injective, each $\lambda_s$ must be injective, and
hence $\lambda_0, \dots, \lambda_{n-1} \in \Sym(\Omega_m)$.

Fix $s \in \Omega_m$, and put $\lambda = \lambda_s$. If the $(i,j)$
entry of $P_{\lambda}$ is nonzero, then $\lambda(j) = i$, so
$\pi ( j n + s ) = \lambda_s(j) n + \rho(s) = i n + \rho(s)$.
Thus the $(i n + \rho(s), j n + s)$ entry of $P_\pi$ is $1$, so the 
$(i n + \rho(s), j n + s)$ entry of $A \dottimes P_\rho$ is nonzero,
and hence $a_{i,j} \ne 0$. Therefore
$\support P_{\lambda} \subseteq \support A$.
\qed\enddemo


\bigskip
\proclaim{4. Proofs of Main Theorems}
\endproclaim

\demo{Proof of {Theorem~1.1}}
Clearly any convex combination of panmagic
permutation matrices is panstochastic.
Conversely, suppose $A \in \gl(5,\reals)$ is panstochastic.
According to
the discussion immediately following {Theorem~2.1},
$
\Pi_5 = \Lambda_5 =
\setof{\pi_{2x+c}}{c \in \Omega_5} \cup
\setof{\pi_{3x+d}}{d \in \Omega_5}
$.
In \cite{12}, Thompson showed that the 
matrices $\setof{P_{\pi}}{\pi \in \Pi_5}$ span $\Pan(5,\reals)$.
Thus there are real numbers $\alpha_c$ and $\beta_d$
such that
$$
\eqalign{
A =& \sum \limits_{c=0}^4 \alpha_c P_{2x+c} +
\sum \limits_{d=0}^4 \beta_d P_{3x+d} \cr
= &
\left [
\matrix
\alpha_0+\beta_0 & \alpha_3+\beta_2 & \alpha_1+\beta_4
& \alpha_4+\beta_1 & \alpha_2+\beta_3 \cr
\alpha_1+\beta_1 & \alpha_4+\beta_3 & \alpha_2+\beta_0
& \alpha_0+\beta_2 & \alpha_3+\beta_4 \cr
\alpha_2+\beta_2 & \alpha_0+\beta_4 & \alpha_3+\beta_1
& \alpha_1+\beta_3 & \alpha_4+\beta_0 \cr
\alpha_3+\beta_3 & \alpha_1+\beta_0 & \alpha_4+\beta_2
& \alpha_2+\beta_4 & \alpha_0+\beta_1 \cr
\alpha_4+\beta_4 & \alpha_2+\beta_1 & \alpha_0+\beta_3
& \alpha_3+\beta_0 & \alpha_1+\beta_2 \cr
\endmatrix 
\right ] . \cr }
\eqno{(1)}
$$
If a cyclic permutation of rows or columns 
is applied 
to a matrix of the form $P_{2x+c}$
($P_{3x+d}$, respectively), the result is 
another matrix
of the form $P_{2x+c}$ ($P_{3x+d}$, respectively).
Hence we can assume without loss of generality that
$\alpha_0 = \min \set{\alpha_c}$ and
$\beta_0 = \min \set{\beta_d}$. Also, 
$\alpha_0+\beta_0 \ge 0$ and
$\sum_c \alpha_c + \sum_d \beta_d = 1$
because $A$ is panstochastic. If both $\alpha_0$ and
$\beta_0$ are nonnegative, then $(1)$ expresses
$A$ as a convex combination of panmagic permutation 
matrices. If $\beta_0 < 0$, then
$$
A = \sum_{c=0}^4 (\alpha_c+\beta_0) P_{2x+c}
+ \sum_{d=0}^4 (\beta_d - \beta_0) P_{3x+d}
$$
is a representation of the required form, and a
similar expression can be obtained if $\alpha_0 < 0$.
\qed\enddemo


We have already dispensed with the case $\gcd(n,6)>1$
in {Theorem~1.2}. The remainder of the proof uses the following
two lemmas.

\proclaim{Lemma~4.1}
If $n>1$ and $\gcd(n,30)=1$, then there is some
panstochastic $n \times n$ matrix that is not a convex 
combination of panmagic permutation matrices.
\endproclaim

\demo{Proof}
The $7 \times 7$ matrix
$$
A = 
{1 \over 2 }
\left [
\matrix
0 & 0 & 0 & 0 & 1 & 1 & 0 \cr
1 & 0 & 1 & 0 & 0 & 0 & 0 \cr
1 & 0 & 0 & 0 & 0 & 1 & 0 \cr
0 & 0 & 0 & 2 & 0 & 0 & 0 \cr
0 & 1 & 0 & 0 & 0 & 0 & 1 \cr
0 & 0 & 0 & 0 & 1 & 0 & 1 \cr
0 & 1 & 1 & 0 & 0 & 0 & 0 \cr
\endmatrix
\right ]
$$
is panstochastic. Suppose $\pi = \pi_{ax+b} \in
\Lambda_7 = \Pi_7$
and
$\support P_\pi \subseteq \support A$.
Then $b = \pi(0) \in \set{1,2}$ and
$\pi(3)=3a+b=3$. If $b=1$, then
$a=3$, so $\pi(2)=0$, which contradicts 
$\support P_\pi \subseteq \support A$.
On the other hand, if $b=2$, then
$a=5$, so $\pi(1)=0$, which again
gives a contradiction.
Thus the assertion of the lemma holds
when $n=7$.

For the rest of the proof we
suppose $n > 7$ and $\gcd(n,30)=1$, so $n \ge 11$.
Observe that 
$\pi_{2x+1}$ and $\pi_{2x-4}$ are affine
panmagic permutations of degree $n$ by {Lemma~2.2}. Define 
$A_0 = (1/2)(P_{2x+1}+P_{2x-4})$, so $A_0$ is panstochastic 
and takes the form
$$
A_0 =
{1 \over 2} 
\left [
\matrix \format \, \r & \quad \r & \quad \r &
\quad \r & \quad \c & \quad \c & \quad \c \, \\
0 & 0 & 1 & 0 & * & \hdots & * \cr
1 & 0 & 0 & 0 & * & \hdots & *\cr
0 & 0 & 0 & 1 & * & \hdots & * \cr
0 & 1 & 0 & 0 & * & \hdots & * \cr
* & * & * & * & * & \hdots & * \cr
\vdots & \vdots & \vdots & \vdots & \vdots & \ddots & \vdots \cr
* & * & * & * & * & \hdots & * \cr
\endmatrix
\right ].
$$
Consider the $n \times n$ matrix
$$
A_1 =
{1 \over 2} 
\left [
\matrix \format \r & \quad \r & \quad \r &
\quad \r & \quad \c & \quad \c & \quad \c \, \\
0 & 1 & -1 & 0 & 0 & \hdots & 0 \cr
-1 & 0 & 0 & 1 & 0 & \hdots & 0 \cr
1 & 0 & 0 & -1 & 0 & \hdots & 0 \cr
0 & -1 & 1 & 0 & 0 & \hdots & 0 \cr
0 & 0 & 0 & 0 & 0 & \hdots & 0 \cr
\vdots & \vdots & \vdots & \vdots & \vdots & \ddots & \vdots \cr
0 & 0 & 0 & 0 & 0 & \hdots & 0 \cr
\endmatrix
\right ],
$$
with zero entries except in rows 
$0$ through $3$ and
columns $0$ through $3$. Observe that
$A_1 \in \Pan(n,\reals)$
and
$\mu(A_1) = 0$.
Define
$A = A_0 + A_1$, so $A$ is
panstochastic and has the form
$$
A =
{1 \over 2} 
\left [
\matrix \format \, \r & \quad \r & \quad \r &
\quad \r & \quad \c & \quad \c & \quad \c \, \\
0 & 1 & 0 & 0 & * & \hdots & * \cr
0 & 0 & 0 & 1 & * & \hdots & *\cr
1 & 0 & 0 & 0 & * & \hdots & * \cr
0 & 0 & 1 & 0 & * & \hdots & * \cr
* & * & * & * & * & \hdots & * \cr
\vdots & \vdots & \vdots & \vdots & \vdots & \ddots & \vdots \cr
* & * & * & * & * & \hdots & * \cr
\endmatrix
\right ].
$$

Assume that $A$ is a convex combination of panmagic permutation
matrices. Then there exists some $\pi \in \Pi_n$ such that
$$
\pi(0)=2
\quad \text{ and }\quad
\support P_\pi \subseteq \support A .
\eqno{(2)}
$$
We show that the existence of such a permutation $\pi$ leads to
a contradiction.

Notice that from {(2)} we have 
$$
\pi(1) \in \left \{ 0, n-2 \right \}, \quad
\pi(2) \in \left \{ 3, 5 \right \}, \quad
\pi(3) \in \left \{1, 7 \right \} .
\eqno{(3)}
$$
Put
$\overline{\Omega}_n=\Omega_n \setminus \left \{ 0,1,2,3 \right \}$.
If $\ell \in \overline{\Omega}_n$, then either 
$\pi(\ell) \equiv 2 \ell + 1 \mymod{n}$
or 
$\pi(\ell) \equiv 2 \ell -4 \mymod{n}$
by {(2)}. Define subsets $J$, $K$ of 
$\overline{\Omega}_n$ as follows:
$$
J = 
\bigsetof{j \in \overline{\Omega}_n}
{ \pi(j) \equiv 2 j + 1 \mymod{n}} ,
$$
$$
K = 
\bigsetof{k \in \overline{\Omega}_n}
{ \pi(k) \equiv 2 k - 4 \mymod{n}} .
$$
Since $\pi$ is panmagic, we have
$$
\pi(i)-i \not\equiv \pi(j)-j
\mymod{n}
\quad\text{whenever}\quad
i , j \in \Omega_n \text{ and } i \ne j
\eqno{(4)}
$$
by {Lemma~2.1}.
Moreover,
$$
\text{if $j \in J$ and 
$j < n-5$, then $j+5 \in J$}.
\eqno{(5)}
$$
Indeed, if $j \in J$ and 
$j+5 \in K$, then 
$$
\eqalign{
\pi(j) - j \, \equiv \, \left ( 2j+1 \right ) - j 
\equiv \, & \left ( 2(j+5) - 4 \right) - (j+5) \cr
\equiv \, & \, \pi(j+5)-(j+5) \mymod{n}, \cr
}
$$
contradicting {(4)}. 
It follows that 
$$
\text{
if $k \in K$ and $k \ge 9$,
then $k-5 \in K$.}
\eqno{(6)}
$$
since $J$ and $K$
form a partition of $ \overline{\Omega}_n$.

The following statements {(7)} -- {(10)} 
all follow from {(4)}.
$$
\eqalignno{
\text{If $n-1 \in J$,} & \ \text{then $4 \in J$.}& {(7)} \cr
\text{If $n-4 \in J$,} & \ \text{then $n-2 \in K$.}& {(8)} \cr
\text{If $n-3 \in J$,} & \ \text{then $8 \in J$.}& {(9)} \cr
\text{If $7 \in K$,} & \ \text{then $5 \in J$.}& {(10)} \cr
}
$$
For {(7)}, suppose $n-1 \in J$ and $4 \in K$. Then
$\pi(n-1) - (n-1) \equiv 0 \equiv \pi(4)-4 \mymod{n}$,
a contradiction. For {(8)}, suppose $n-2 \in J$ and
$n-4 \in J$. From the definition of $J$ we have
$\pi(n-2)-(n-2) \equiv n-1 \mymod{n}$ and
$\pi(n-4)-(n-4) \equiv n-3 \mymod{n}$.
However, by {(3)},
$\pi(1)-1 \equiv n-1 \mymod{n}$ or
$\pi(1)-1 \equiv n-3 \mymod{n}$, which is a contradiction.
The proofs of {(9)} and {(10)} are similar to that
of {(8)}, using the values of $\pi(3)$ and $\pi(2)$,
respectively. 

We have 
$$
j \in J
\quad
\text{whenever}
\quad
4 \le j < n
\text{ and }
j \equiv 1 \mymod{5}.
\eqno{(11)}
$$
Indeed, by {(5)}, it is enough to show $6 \in J$.
For this, observe $\pi(0) - 0 = 2$, so $\pi(6) - 6 \ne 2$
by {(4)}, and hence $6 \not \in K$.

In addition, we have
$$
{n+1 \over 2} \in K.
\eqno{(12)}
$$
Indeed, if $(n+1)/2 \in J$, then 
$\pi \left( (n+1)/2 \right) \equiv n+2 \equiv 2 \equiv \pi(0) \mymod{n}$,
a contradiction.

We finish the argument by cases, according to the congruence class
of $n$ modulo $5$.

\noindent{\it Case 1:}
$n \equiv 1 \mymod{5}$.\ \ 
In this case $(n+1)/2 \equiv 1 \mymod{5}$,
and thus $(n+1)/2 \in J$ by {(11)}.
This contradicts {(12)}.

\noindent{\it Case 2:}
$n \equiv 2 \mymod{5}$.\ \ 
In this case, $n-1 \equiv 1 \mymod{5}$,
and thus $n-1 \in J$ by {(11)}.
By {(7)}, $4 \in J$, and by
{(5)}, $j \in J$ whenever
$4 \le j < n$ and $j \equiv 4 \mymod{5}$.
Since $(n+1)/2 \equiv 4 \mymod{5}$,
$(n+1)/2 \in J$, contradicting {(12)}.

\noindent{\it Case 3:}
$n \equiv 3 \mymod{5}$.\ \ 
Since $(n+1)/2 \equiv 2 \mymod{5}$, it follows from
{(12)} and {(6)} that $7 \in K$, and
thus $\pi(7)\equiv 10\mymod{n}$. However, if $j = (n+9)/2$,
then $4 \le j < n$ and $j \equiv 1 \mymod{5}$,
so $j \in J$ by {(11)}, and hence
$\pi(j) \equiv 2j+1 \equiv 10 \mymod{n}$,
which is a contradiction. 

\noindent{\it Case 4:}
$n \equiv 4 \mymod{5}$.\ \ 
In this case, $n-3 \equiv 1 \mymod{5}$, and so
$n-3 \in J$ by {(11)}. Hence 
$8 \in J$ by {(9)}, and hence 
$n-1 \in J$ by {(5)} since
$n-1 \equiv 8 \mymod{5}$. Thus
$4 \in J$ by {(7)}, and so
$j \in J$ whenever $4 \le j < n$ and
$j \equiv 4 \mymod{5}$ by {(5)}.
Since $n-4 \equiv 4 \mymod{5}$,
we conclude $n-4 \in J$. By
{(8)}, $n-2 \in K$, and since
$n-2 \equiv 2 \mymod{5}$, we have
$7 \in K$ by {(6)}. By
{(10)}, $5 \in J$, and by
{(5)}, $j \in J$ whenever 
$4 \le j < n$ and $j \equiv 0 \mymod{5}$.
Since $(n+1)/2 \equiv 0 \mymod{5}$,
we have $(n+1)/2 \in J$,
contradicting {(12)}.

We have arrived at a contradiction in each case.
\qed\enddemo


\proclaim{Lemma~4.2}
Suppose there is some $m \times m$ panstochastic matrix that is
not a convex combination of panmagic permutation matrices
and that $n$ is a positive integer with $\gcd(n,6)=1$. Then there is
some $(m n) \times (m n)$ panstochastic matrix that is not a convex
combination of panmagic permutation matrices.
\endproclaim

\demo{Proof}
There exists a panmagic permutation $\rho$ of degree $n$, say
$\rho = \pi_{2x}$. Among all panstochastic $m \times m$
matrices that are not convex combinations of panmagic
permutation matrices, choose one, say $A$, with the maximum 
number of zero entries. The matrix $B = A \dottimes P_\rho$ is 
then a panstochastic $(m n) \times (m n)$ matrix by
{Corollary~3.1}. Assume that $B$ is a convex combination
of panmagic permutation matrices. Then there is some
panmagic permutation $\pi$ of degree $m n$ such that
$\support P_\pi \subseteq \support B$. By {Theorem~3.3}, there are 
$\lambda_0, \dots \lambda_{n-1} \in \Sym(\Omega_m)$ such that
$
\pi = (\lambda_0, \dots, \lambda_{n-1}) \wr \rho
$
and
$\support P_{\lambda_s} \subseteq \support A$ for all $s \in \Omega_m$.
Then $\lambda_0, \dots, \lambda_{n-1}$ are panmagic by
{Theorem~3.2}. Fix some $s \in \Omega_m$ and put
$\lambda = \lambda_s$, so $\support P_\lambda \subseteq \support A$. 
Define
$$
a = \min \, \setof{a_{i,j}}{\lambda(j)=i}.
$$
We cannot have $a=1$, for otherwise $A = P_{\lambda}$ is a panmagic
permutation matrix. Therefore $0 < a < 1$, and so the matrix
$$
C = {1 \over 1 - a} A \, - \, { a \over 1-a } P_{\lambda} 
$$
is a panstochastic $m \times m$ matrix with a greater number of 
zero entries than $A$. By the choice of $A$, $C$ is a convex 
combination of panmagic permutation matrices, and it follows that
$
A = a P_\lambda + (1 - a) C
$
is also a convex combination of panmagic permutation matrices,
so we have a contradiction. Therefore $B$ is not a convex 
combination of panmagic permutation matrices.
\qed\enddemo


\demo{Proof of {Theorem~1.2}}
By {Lemma~4.1} and {Lemma~4.2}, it suffices to prove that there
is some $25 \times 25$ panstochastic matrix that is not a 
convex combination of panmagic permutation matrices.
The matrix 
%
%
\def\sz{\text{\smallfont 0}}
\def\so{\text{\smallfont 1}}
\def\ssp{{\hskip 4 pt}}
$$
\spreadmatrixlines{-3pt}
A =
{1 \over 2}
\left [
\matrix \format
\c \ssp & \c \ssp & \c \ssp & \c \ssp & \c \ssp & \c \ssp & \c \ssp & \c \ssp & \c \ssp & 
\c \ssp & \c \ssp & \c \ssp & \c \ssp & \c \ssp & \c \ssp & \c \ssp & \c \ssp & \c \ssp & 
\c \ssp & \c \ssp & \c \ssp & \c \ssp & \c \ssp & \c \ssp & \c \\
\sz&\so&\sz&\sz&\sz&\sz&\sz&\sz&\sz&\sz&\sz&\sz&\so&\sz&\sz&\sz&\sz&
\sz&\sz&\sz&\sz&\sz&\sz&\sz&\sz\cr
\sz&\sz&\sz&\so&\sz&\sz&\sz&\sz&\sz&\sz&\sz&\so&\sz&\sz&\sz&\sz&\sz&
\sz&\sz&\sz&\sz&\sz&\sz&\sz&\sz\cr
\so&\sz&\sz&\sz&\sz&\sz&\sz&\sz&\sz&\sz&\sz&\sz&\sz&\so&\sz&\sz&\sz&
\sz&\sz&\sz&\sz&\sz&\sz&\sz&\sz\cr
\sz&\sz&\so&\sz&\sz&\sz&\sz&\sz&\sz&\sz&\sz&\sz&\sz&\sz&\sz&\sz&\sz&
\so&\sz&\sz&\sz&\sz&\sz&\sz&\sz\cr
\sz&\sz&\sz&\sz&\sz&\sz&\sz&\sz&\sz&\sz&\sz&\sz&\sz&\sz&\so&\sz&\sz&
\sz&\sz&\sz&\sz&\so&\sz&\sz&\sz\cr
\sz&\sz&\so&\sz&\sz&\sz&\sz&\sz&\sz&\sz&\sz&\sz&\sz&\sz&\sz&\sz&\so&
\sz&\sz&\sz&\sz&\sz&\sz&\sz&\sz\cr
\sz&\sz&\sz&\sz&\sz&\so&\sz&\sz&\sz&\sz&\sz&\sz&\sz&\sz&\sz&\so&\sz&
\sz&\sz&\sz&\sz&\sz&\sz&\sz&\sz\cr
\sz&\sz&\sz&\so&\sz&\sz&\sz&\sz&\sz&\sz&\sz&\sz&\sz&\sz&\sz&\sz&\sz&
\sz&\sz&\so&\sz&\sz&\sz&\sz&\sz\cr
\sz&\sz&\sz&\sz&\sz&\sz&\sz&\sz&\sz&\sz&\sz&\sz&\sz&\sz&\sz&\so&\so&
\sz&\sz&\sz&\sz&\sz&\sz&\sz&\sz\cr
\sz&\sz&\sz&\sz&\so&\sz&\sz&\so&\sz&\sz&\sz&\sz&\sz&\sz&\sz&\sz&\sz&
\sz&\sz&\sz&\sz&\sz&\sz&\sz&\sz\cr
\sz&\sz&\sz&\sz&\sz&\sz&\sz&\sz&\sz&\sz&\sz&\sz&\sz&\sz&\sz&\sz&\sz&
\so&\sz&\sz&\sz&\sz&\sz&\so&\sz\cr
\sz&\sz&\sz&\sz&\sz&\so&\sz&\sz&\so&\sz&\sz&\sz&\sz&\sz&\sz&\sz&\sz&
\sz&\sz&\sz&\sz&\sz&\sz&\sz&\sz\cr
\sz&\sz&\sz&\sz&\sz&\sz&\so&\sz&\sz&\sz&\sz&\sz&\sz&\sz&\sz&\sz&\sz&
\sz&\so&\sz&\sz&\sz&\sz&\sz&\sz\cr
\sz&\sz&\sz&\sz&\sz&\sz&\so&\sz&\sz&\so&\sz&\sz&\sz&\sz&\sz&\sz&\sz&
\sz&\sz&\sz&\sz&\sz&\sz&\sz&\sz\cr
\sz&\sz&\sz&\sz&\sz&\sz&\sz&\sz&\sz&\sz&\sz&\sz&\sz&\sz&\sz&\sz&\sz&
\sz&\so&\so&\sz&\sz&\sz&\sz&\sz\cr
\sz&\sz&\sz&\sz&\sz&\sz&\sz&\so&\sz&\sz&\sz&\sz&\sz&\sz&\sz&\sz&\sz&
\sz&\sz&\sz&\sz&\sz&\sz&\sz&\so\cr
\sz&\sz&\sz&\sz&\sz&\sz&\sz&\sz&\sz&\sz&\sz&\sz&\sz&\sz&\sz&\sz&\sz&
\sz&\sz&\sz&\so&\sz&\so&\sz&\sz\cr
\so&\sz&\sz&\sz&\sz&\sz&\sz&\sz&\so&\sz&\sz&\sz&\sz&\sz&\sz&\sz&\sz&
\sz&\sz&\sz&\sz&\sz&\sz&\sz&\sz\cr
\sz&\sz&\sz&\sz&\sz&\sz&\sz&\sz&\sz&\sz&\sz&\sz&\sz&\so&\sz&\sz&\sz&
\sz&\sz&\sz&\sz&\so&\sz&\sz&\sz\cr
\sz&\sz&\sz&\sz&\sz&\sz&\sz&\sz&\sz&\so&\so&\sz&\sz&\sz&\sz&\sz&\sz&
\sz&\sz&\sz&\sz&\sz&\sz&\sz&\sz\cr
\sz&\sz&\sz&\sz&\sz&\sz&\sz&\sz&\sz&\sz&\sz&\sz&\sz&\sz&\sz&\sz&\sz&
\sz&\sz&\sz&\so&\sz&\so&\sz&\sz\cr
\sz&\sz&\sz&\sz&\sz&\sz&\sz&\sz&\sz&\sz&\so&\sz&\sz&\sz&\so&\sz&\sz&
\sz&\sz&\sz&\sz&\sz&\sz&\sz&\sz\cr
\sz&\sz&\sz&\sz&\sz&\sz&\sz&\sz&\sz&\sz&\sz&\sz&\so&\sz&\sz&\sz&\sz&
\sz&\sz&\sz&\sz&\sz&\sz&\so&\sz\cr
\sz&\so&\sz&\sz&\sz&\sz&\sz&\sz&\sz&\sz&\sz&\so&\sz&\sz&\sz&\sz&\sz&
\sz&\sz&\sz&\sz&\sz&\sz&\sz&\sz\cr
\sz&\sz&\sz&\sz&\so&\sz&\sz&\sz&\sz&\sz&\sz&\sz&\sz&\sz&\sz&\sz&\sz&
\sz&\sz&\sz&\sz&\sz&\sz&\sz&\so\cr
\endmatrix 
\right ] 
$$
is panstochastic; it can be obtained by averaging the matrices for
$\pi_{2x+1}$ and a non-affine panmagic permutation and then
adjusting the entries in rows 0 through 3 and columns 0 through 3,
as in the proof of {Lemma~4.1}. Suppose $A$ is a convex 
combination of panmagic permutation matrices. Thus there is some
panmagic permutation $\pi$ on $\Omega_{25}$ such that
$\pi(0) = 2$ and $\support P_\pi \subseteq \support A$.
Since $\pi(0)=2$, we cannot have $\pi(13)=2$, and hence $\pi(13)=18$.
Thus $\pi(21) \ne 18$, so $\pi(21) = 4$. Also, $\pi(0) \ne 17$,
so $\pi(8) = 17$. Therefore
$\pi(8)+8 \equiv 0 \equiv \pi(21)+21 \mymod{25}$,
so $\pi$ is not panmagic by {Lemma~2.1}, and a contradiction is
reached. 
\qed\enddemo


\bigskip
\noindent
{\it Indiana University South Bend, 
South Bend, IN 46634}
\hfil\break\noindent
{\it dalvis\@iusb.edu}, {\it mkinyon\@iusb.edu}


\enddocument
\bye